\documentclass{amsart}
\usepackage[all]{xy}
\usepackage{lipsum}
\usepackage{amsfonts}
\usepackage{graphicx}
\usepackage{amsmath}
\usepackage{amssymb}
\usepackage[english]{babel}
\usepackage{tikz-cd}
\usepackage{hyperref}

\newtheorem{theorem}{Theorem}
\newtheorem{lemma}{Lemma}
\newtheorem{proposition}{Proposition}
\newtheorem{corollary}{Corollary}

\newtheorem{remark}[theorem]{Remark}

\begin{document}

\title{Classes of Preradicals and Relative Injectivity}

\author{Frank Murphy-Hernandez}
\address{Facultad de Ciencias, UNAM, Mexico City}
\email{murphy@ciencias.unam.mx}

\author{José Ríos}
\address{Instituto de Matemáticas, UNAM, Mexico City}
\email{ jrios@matem.unam.mx }


\date{\today}


\begin{abstract}
This paper analyzes new classes of preradicals defined as weak forms of left exact and idempotent preradicals. We introduce prehereditary preradicals to generalize the hereditary property, and essentially idempotent and weakly idempotent preradicals to generalize idempotency. Our motivation is to study the concept of relative injectivity. We show that these new classes facilitate the extension of classical results on injectivity with respect to a hereditary torsion theory to a broader context, requiring weaker assumptions.
\end{abstract}

\maketitle

\section{Introduction}
In this work, we introduce new classes of preradicals and study their relationship with injectivity relative to a preradical.

We first present the class of \emph{essentially idempotent preradicals} as a generalization of idempotent preradicals. As expected, these preradicals preserve several properties of their idempotent counterparts; for example, the class of essentially idempotent preradicals is closed under suprema, just like the class of idempotent preradicals, and their pretorsion-free classes are closed under extensions.

Subsequently, we introduce a generalization of essentially idempotent preradicals: \emph{weakly idempotent preradicals}. This class also preserves some of the properties of the class of idempotent preradicals that the essentially idempotent ones share. This class is introduced to naturally define its complement: the class of \emph{strongly nilpotent preradicals}.

The class of \emph{prehereditary preradicals} is introduced in this article as a generalization of left exact or hereditary preradicals. We prove that there exists a correspondence between linear filters and prehereditary preradicals. This correspondence is not a lattice isomorphism but rather an epimorphism from the class of prehereditary preradicals to linear filters.

Preradicals provide a suitable context for defining essentiality relative to a preradical, which preserves almost all properties of usual essentiality. We study the concept of purity relative to a preradical and obtain analogs of the classical results that hold when the preradical is a left exact radical. We observe that with an idempotent radical, all usual properties are recovered.

We investigate the concept of injectivity relative to a preradical and find that for good properties, it is sufficient for the preradical to be an idempotent radical. However, if one desires criteria analogous to Baer's Criterion, the property of being prehereditary plays a crucial role.

It is always possible to define the relative injective hull with respect to a preradical, but when the preradical is an idempotent radical, it is unique with respect to certain properties. We also define pseudocomplemented submodules relative to a preradical, which provide conditions for determining when two preradicals have the same class of injective modules.

We study torsion-free injective modules, called absolutely pure modules, and continue with the next class of preradicals: \emph{autocostable preradicals}, which are those whose pretorsion-free classes are closed under relative injective hulls. We show that under certain requirements, this property implies costability.

Finally, we define an assignment for any left module with respect to a preradical that, when the preradical is a left exact radical, coincides with the localization functor. Certain properties of this assignment provide information about the preradical even when it is not a left exact radical.

\section{Preliminaries}

In this paper $R$ always will denote an associative ring with $1$, and by $R$-Mod it will understand the category of unitary left modules over $R$, an excelent reference of the category of left $R$-modules is (10). Let $M$ be an $R$-module and $N$ a submodule of $M$, if $N$ is essential in $M$ this is denoted by $N\unlhd M$. The injective hull of $M$ is denoted by $E(M)$.

A preradical $\sigma$ over $R$-Mod is a subfunctor of the identity, references for preradicals are (1) and (9). The class of the preradicals is ordered punctually, that is, if $\sigma$ and $\tau$ are preradicals over $R$-Mod, $\sigma\leq\tau$ if $\sigma(M)\leq\tau(M)$ for any left $R$-module $M$, with this order the class of preradicals becomes a big complete lattice, it also, has to operations the product and the coproduct, if a preradical is idempotent under de product it is called an idempotent preradical and if is idempotent under de coproduct it is called a radical. The class of all idempotent preradicals over $R$-Mod will be denote by $R$-id and the class of all radicals over $R$-Mod will be denote by $R$-rad. A class of left $R$-modules is called a precursion class if it is closed under coproducts and quotients and it is called a pretorsion free class if it is closed under products and submodules. If $\sigma$ is a preradical the class ${\mathbb{T}}_{\sigma}=\{M\mid\sigma(M)=M\}$ is a pretorsion class and if ${\mathbb{T}}$ is a pretorsion class then the assignation $\sigma_{{\mathbb{T}}}(M)=\sum\{N\leq M\mid N\in{\mathbb{T}}\}$ for any left $R$-module $M$ is an idempotent preradical, this is a bijective correspondence between the pretorsion classes and idempotent preradicals. In the same way, if $\sigma$ is preradical then ${\mathbb{F}}_{\sigma}=\{M\mid\sigma(M)=0\}$ is a pretorsion free class and if ${\mathbb{F}}$ is pretorsion free class then assignation $\sigma^{{\mathbb{F}}}(M)=\bigcap\{N\leq M\mid\sigma(M/N)=0\}$ for any left $R$-module $M$ is a radical, this is a bijective correspondence between the pretorsion free classes and the radicals. The elements of ${\mathbb{T}}_{\sigma}$ are called $\sigma$-torsion modules and the elements of ${\mathbb{F}}_{\sigma}$ are called $\sigma$-torsion free modules. A class of modules ${\mathcal{C}}$ is closed under extensions if for any short exact sequence $0\longrightarrow M^{\prime}\longrightarrow M\longrightarrow M^{\prime\prime} \longrightarrow 0$ with $M^{\prime},M^{\prime\prime}\in{\mathcal{C}}$ then $M\in{\mathcal{C}}$. A pretorsion class which is closed under extensions is called a torsion class and a pretorsion free class which is closed under extensions is called a torsion free class.If $\sigma$ is an idempotent preradical then ${\mathbb{T}}_{\sigma}$ is a torsion class and If $\sigma$ is a radical then ${\mathbb{F}}_{\sigma}$ is a torsion free class.

The class of idempotent preradicals is closed under supremum so it is possible for any preradical $\sigma$ to obtain the greatest idempotent preradical below $\sigma$, it results to be $\sigma_{{\mathbb{T}}_{\sigma}}$ and it is usually denoted by $\widehat{\sigma}$. In the same manner the class of radicals is closed under infimum so for any preradical $\sigma$ there is the least radical above $\sigma$, it results to be $\sigma^{{\mathbb{F}}_{\sigma}}$ and it is usually denoted by $\bar{\sigma}$

Let $M$ be a left $R$-module and let $N$ be a submodule of $M$, $N$ is called a fully invariant submodule of $M$ if $f(N)\subseteq N$ for any endomorphism $f$ of $M$. Let $\sigma$ be a preradical, then $\sigma(M)$ is a fully invariant submodule of $M$, moreover, $N$ is a fully invariant submodule of $M$ if and only if there is a preradical $\sigma$ such that $\sigma(M)=N$. If $N$ is a fully invariant submodule of $M$, it is defined:

\[\alpha_{N}^{M}(K)=\sum\{f(N)\mid f:M\longrightarrow K\}\]

\[\omega_{N}^{M}(K)=\bigcap\{f^{-1}(N)\mid f:K\longrightarrow M\}\]

for any left $R$-module $K$. It is easy to see that $\alpha_{N}^{M}(M)=N$ and $\omega_{N}^{M}(M)=N$, moreover, for any preradical $\sigma$, $\sigma(M)=N$ if and only if $\alpha_{N}^{M}\leq\sigma\leq\omega_{N}^{M}$. If $N$ is a fully invariant submodule of $M$, is defined $\widehat{N}=\overline{\alpha_{N}^{M}(M)}$, then $\alpha_{N}^{M}$ is an idempotent preradical if and only if $N=\widehat{N}$, in the same way it is defined $\bar{N}=\omega_{N}^{M}(M)$ then $\omega_{N}^{M}$ is a radical if and only if $N=\bar{N}$. Let $S$ be a simple left $R$-module, then $S$ is injective if and only is $\alpha_{S}^{E(S)}$ is idempotent.. The standard references for lattice aspects of preradicals and the alphas and omegas preradicals are the papers (6), (7) and (8).

A class of modules is called hereditary if it is closed under submodules, a preradical $\sigma$ is called hereditary if it is idempotent and ${\mathbb{T}}_{\sigma}$ is hereditary. A preradical is hereditary if and only if it is left exact. The class of all left exact preradicals over $R$-Mod will be denoted by $R$-lep, this class is closed under infimum so for any preradical $\sigma$ there is the least left exact preradical $\widetilde{\sigma}$ above $\sigma$, $\widetilde{\sigma}$ is easy to describe, $\widetilde{\sigma}(M)=\sigma(E(M))\cap M$ for any left $R$-module $M$. A set of left ideals ${\mathbb{I}}$ that satifies:

* *If $I\in{\mathbb{I}}$ and $I\subseteq J\leq R$ then $J\in{\mathbb{I}}$.
* *If $I,J\in{\mathbb{I}}$ then $I\cap J\in{\mathbb{I}}$.
* *If $a\in R$ and $I\in{\mathbb{I}}$ then $(I:a)\in{\mathbb{I}}$.

is called a left linear filter. If $\sigma$ be a left exact preradical then it is defined ${\mathbb{I}}_{\sigma}=\{I\leq R\mid\sigma(R/I)=R/I\}$ and if ${\mathbb{I}}$ is a linear filter then it is defined $\sigma_{{\mathbb{I}}}(M)=\{x\in M\mid ann(x)\in{\mathbb{I}}\}$, this is a bijective correspondence between the left exact preradicals and the left linear filters. A left linear filter ${\mathbb{I}}$ is called a left Gabriel filter if it satisfies: If $I\in{\mathbb{I}}$ and $J\leq R$ is such that for any $a\in I$$(J:a)\in{\mathbb{I}}$ then $J\in{\mathbb{I}}$. The previous correspondence induces a bijective correspondence between the left exact radicals and the left Gabriel filters.

Let $E$ and $E^{\prime}$ be injective left $R$-modules, it is said that $E$ and $E^{\prime}$ are related if there is an imbedding of $E$ in a product of copies of $E^{\prime}$ and there is an embedding of $E^{\prime}$ in a product of copies of $E$, it is easy to see that this an equivalence relation, an a class of equivalence is a called an hereditary torsion theory, a good reference is (5). There is a bijective correspondence between hereditary torsion theories and the left exact radicals.

\section{Essentially Idempotent Preradicals}
Let $\sigma$ be a preradical over $R$-Mod, $\sigma$ is called essentially idempotent if $\sigma(M)\neq 0$ implies $\widehat{\sigma}(M)\neq 0$ for any left $R$-module $M$. It is observed that any idempotent preradical is essentially idempotent so the property of being essentially idempotent is a generalization of being idempotent. The class of all the essentially idempotent preradicals over $R$-Mod is denoted by $R$-eid , the last remark could be restated as $R$-id$\subseteq R$-eid. From this last fact it could happens that $R$-eid is not a set, since $R$-id is not always a set. As the supremum of a family of idempotent preradicals is idempotent, it is expected that the same happens for essentially idempotent preradicals.

\begin{proposition}
Let $\{\sigma_{i}\}_{i\in I}$ be a family of essentially idempotent preradicals over $R$-Mod. Then $\bigvee_{i\in I}\sigma_{i}$ is an essentially idempotent preradical.
\end{proposition}

\textbf{Proof.} Let $M$ be a left $R$-module with $(\bigvee_{i\in I}\sigma_{i})(M)\neq 0$, then there is $i\in I$ with $\sigma_{j}(M)\neq 0$ and by hypothesis it follows that $\widehat{\sigma}_{j}(M)\neq 0$, so $(\bigvee_{i\in I}\widehat{\sigma}_{i})(M)\neq 0$ and since $\bigvee_{i\in I}\widehat{\sigma}_{i}\leq\widehat{\bigvee_{i\in I}\sigma_{i}}$, $\widehat{\bigvee_{i\in I}\sigma_{i}}(M)\neq 0$ as desired. $\blacksquare$

By the last proposition for any preradical $\sigma$ over $R$-Mod, it is possible to construct an essentially idempotent preradical $\sigma^{\circ}$ as the supremum of all essentially idempotent preradicals below $\sigma$, by the previous proposition $\sigma^{\circ}$ is an essentially idempotent preradical. In fact $\sigma^{\circ}$ is the greatest essentially idempotent below $\sigma$. It is observed that a preradical $\sigma$ is essentially idempotent if and only if $\sigma^{\circ}=\sigma$. Also it is important to remember that the class $R$-id is no closed under infimum, even finite ones, this pathology is preserved by the class $R$-eid, consider the next example. Let $R$ be the ring of the integers, let $\sigma$ be the socle and $\tau$ the divisible part, as $\sigma$ and $\tau$ are idempotents they are essentially idempotents, but $(\sigma\wedge\tau)(\mathbb{Z}_{p^{\infty}})=\mathbb{Z}_{p}$ for any prime $p$, then $(\sigma\wedge\tau)^{2}(\mathbb{Z}_{p^{\infty}})=0$ which implies that $\widehat{(\sigma\wedge\tau)}(\mathbb{Z}_{p^{\infty}})=0$, so $\sigma\wedge\tau$ is not essentially idempotent, from here it is seen that the infimum of idempotent preradicals is not essentially idempotent, which also implies that the infimum of essentially idempotent preradicals is not essentially idempotent. The last said that $R$-eid is not a sublattice of $R$-pr and $R$-id is not a sublattice of $R$-eid.

As always $R$-eid has a natural way to be described as a complete lattice, that is for any family $\{\sigma_{i}\}_{i\in I}$ of essentially idempotent preradicals the supremum is the usual supremum in $R$-pr, but the infimum results $(\bigwedge_{i\in I}\sigma_{i})^{\circ}$. The next proposition tells that the operator ${}^{\circ}$ over $R$-pr is an interior operator.

\begin{proposition}
The assignation ${}^{\circ}:R\text{-pr}\longrightarrow R\text{-pr}$ given by $\sigma\mapsto\sigma^{\circ}$ for any preradical over $R$-Mod $\sigma$ is a monotone, deflatory and idempotent operator over $R$-pr.
\end{proposition}

From the fact that $R$-id$\subseteq R$-eid follows that $\widehat{\sigma}(M)\leq\sigma^{\circ}(M)$ for any left $R$-module $M$. Then

\begin{remark}
Let $\sigma$ be an essentially idempotent preradical over $R$-Mod and let $M$ be a left $R$-module. If $\sigma(M)=M$ then $\sigma^{\circ}(M)=M$.
\end{remark}

\begin{remark}
Let $\sigma$ be a preradical over $R$-Mod. Then $\mathbb{T}_{\sigma}=\mathbb{T}_{\sigma^{\circ}}=\mathbb{T}_{\widehat{\sigma}}$.
\end{remark}

\begin{remark}
Let $\sigma$ be a preradical over $R$-Mod. Then $\widehat{\sigma^{\circ}}=\widehat{\sigma}$.
\end{remark}

It is well known that for any idempotent preradical $\sigma$, the associated radical $\bar{\sigma}$ is an idempotent radical. The next result has the same spirit that this one.

\begin{proposition}
Let $\sigma$ be an essentially idempotent preradical over $R$-Mod. Then $\bar{\sigma}$ is an essentially idempotent radical.
\end{proposition}

\textbf{Proof.} Let $M$ be a left $R$-module with $\bar{\sigma}(M)\neq 0$, then $\sigma(M)\neq 0$ it follows that $\widehat{\sigma}(M)\neq 0$ and as $\widehat{\sigma}(M)\leq(\widehat{\bar{\sigma}})(M)$ the desired result is obtained.

\begin{proposition}
Let $S$ be a simple left $R$-module. If $\alpha_{S}^{E(S)}$ is essentially idempotent then $\alpha_{S}^{E(S)}$ is idempotent.
\end{proposition}

\textbf{Proof.} First it is remembered that $\alpha_{S}^{E(S)}$ is an atom in the lattice $R$-pr, by this $\widehat{\alpha_{S}^{E(S)}}$ has two options to be $\alpha_{S}^{E(S)}$ or $0$, but $\alpha_{S}^{E(S)}(E(S))=S\neq 0$ which means that $\widehat{\alpha_{S}^{E(S)}}(E(S))\neq 0$, so $\widehat{\alpha_{S}^{E(S)}}=\alpha_{S}^{E(S)}$. $\blacksquare$

\begin{corollary}
Let $R$ be a ring. Then $R$ is a V-ring if and only if every atom in $R$-pr is essentially idempotent.
\end{corollary}

\begin{proposition}
Let $\{\sigma_{i}\}_{i\in I}$ be a family of preradicals over $R$-Mod. Then $\widehat{\bigwedge_{i\in I}\sigma_{i}}=\widehat{\bigwedge_{i\in I}\widehat{ \sigma}_{i}}$.
\end{proposition}

\textbf{Proof.} It is observed that $\mathbb{T}_{\bigwedge_{i\in I}\sigma_{i}}=\bigcap_{i\in I}\mathbb{T}_{\sigma_{i} }=\bigcap_{i\in I}\mathbb{T}_{\widehat{\sigma}_{i}}=\mathbb{T}_{\bigwedge_{i\in I }\widehat{\sigma}_{i}}$.

\begin{corollary}
Let $\{\sigma_{i}\}_{i\in I}$ be a family of preradicals over $R$-Mod such that $\bigwedge_{i\in I}\sigma$ is essentially idempotent. Then $\bigwedge_{i\in I}\widehat{\sigma}$ is essentially idempotent.
\end{corollary}

\begin{remark}
Let $\sigma$ be a preradical over $R$-Mod. Then $\sigma$ is an essentially idempotent preradical if and only if $\mathbb{F}_{\widehat{\sigma}}=\mathbb{F}_{\sigma}$.
\end{remark}

From the classical theory of preradical it is a well known fact that when $\sigma$ is an idempotent preradical $\mathbb{F}_{\sigma}$ is closed under extensions. The next remark generalizes the previous fact

\begin{remark}
Let $\sigma$ be an essentially idempotent preradical over $R$-Mod. Then ${\mathbb{F}}_{\sigma}$ is closed under extensions.
\end{remark}

\begin{proposition}
Let $\sigma$ be an essentially idempotent radical over $R$-Mod. Then $\sigma$ is an idempotent radical.
\end{proposition}

\textbf{Proof.} Let $M$ be a left $R$-module. Since $\sigma$ is a radical, also $\sigma^{2}$ is radical, then $\sigma^{2}(M/\sigma^{2}(M))=0$, this implies that $\sigma(M/\sigma^{2}(M))=0$, but $\sigma(M/\sigma^{2}(M))=(\sigma^{2}:\sigma)(M)/\sigma^{2}(M)$, from here it is obtained that $\sigma(M)\leq(\sigma^{2}:\sigma)(M)=\sigma^{2}(M)$ and the desired result is followed.

It is known that for any left $R$-module $M$ and any fully invariant submodule $N$, the preradical $\alpha_{N}^{M}$ is idempotent if and only if $\widehat{N}=N$. In the same spirit it stated the next remark.

\begin{remark}
Let $M$ be a left $R$-module and $N$ a no zero fully invariant submodule of $M$. If $\alpha_{N}^{M}$ is an essentially idempotent preradical then $\widehat{N}\neq 0$.
\end{remark}

It is considered the ring $R=\mathbb{Z}_{4}\times\mathbb{Z}_{4}$ and the ideal $I=\mathbb{Z}_{4}\times 2\mathbb{Z}_{4}$, it is observed that $\alpha_{I}^{R}(0\times\mathbb{Z}_{4})=0\times 2\mathbb{Z}_{4}$ and $\overline{\alpha_{I}^{\widehat{R}}}(0\times\mathbb{Z}_{4})=0$ which means that $\alpha_{I}^{R}$ is not an essentially idempotent preradical, also that$\widehat{I}=\mathbb{Z}_{4}\times 0$. This tells that in general $M$ is not a test module for the preradical $\alpha_{N}^{M}$ to be essentially idempotent.

It is possible to define the dual concept as, a preradical $\sigma$ is essentially coidempotent if $\bar{\sigma}(M)=M$ then $\sigma(M)=M$ for any left $R$-module $M$. All the previous results in their dual versions are preserved.	
\section{Weakly Idempotent Preradicals}
Let $\sigma$ be a preradical over $R$-Mod, $\sigma$ is called weakly idempotent if $\sigma\neq 0$ implies $\widehat{\sigma}\neq 0$. It is stated this way so the $0$ preradical is weakly idempotent. It is noted that a preradical $\sigma$ is weakly idempotent if and if there exists a no zero idempotent preradical $\tau$ such that $\tau\leq\sigma$. Now, let $R$ be the ring of the integers and let $J$ be the Jacobson radical over $R$-Mod, it is well known that $J$ is not idempotent and that $\widehat{J}$ is $d$ the divisible part. Moreover $J(\mathbb{Z}_{p^{k}})=\mathbb{Z}_{p^{k-1}}$ and $d(\mathbb{Z}_{p^{k}})=0$ for any prime $p$ and any integer $k\geq 1$, which means that $J$ is neither essentially idempotent.

It is noted that any essentially idempotent preradical is weakly idempotent then the property of being weakly idempotent is a generalization of being essentially idempotent. The class of all the weakly idempotent preradicals over $R$-Mod is denoted by $R$-wid, the last remark could be restated as $R$-$\operatorname{cid}_{\subseteq}^{-}R$-wid. As happens with $R$-$\operatorname{cid}$, $R$-wid is not necessarily a set, but as $R$-$\operatorname{cid}$, $R$-wid is closed under supremum.

\begin{proposition}
Let $\{\sigma_{i}\}_{i\in I}$ be a family of weakly idempotent preradicals over $R$-Mod. Then $\bigvee_{i\in I}\sigma_{i}$ is an weakly idempotent preradical.
\end{proposition}

By the last proposition for any preradical $\sigma$ over $R$-Mod, it is possible to construct an weakly idempotent preradical $\sigma^{*}$ as the supremum of all weakly idempotent preradicals below $\sigma$, by the previous proposition $\sigma^{*}$ is an essentially idempotent preradical. In fact $\sigma^{*}$ is the greatest weakly idempotent below $\sigma$. It is observed that a preradical $\sigma$ is essentially idempotent if and only if $\sigma^{*}=\sigma$. Also it is important to remember that the class $R$-$\operatorname{cid}$ is closed under infimum, even finite ones, this pathology is preserved by the class $R$-wid, consider the next example. Let $R$ be the ring of the integers, let $\sigma$ be the socle and $\tau$ the divisible part, as $\sigma$ and $\tau$ are no zero idempotents they are essentially idempotents and $\sigma\wedge\tau\neq 0$, but $(\sigma\wedge\tau)(\mathbb{Z}_{p^{\infty}})=\mathbb{Z}_{p}$ for any prime $p$, then $\widetilde{\sigma\wedge\tau}(\mathbb{Z}_{p^{\infty}})=0$, which implies that $\widehat{(\sigma\wedge\tau)}=0$, so $\sigma\wedge\tau$ is not weakly idempotent, from here it is seen that the infimum of idempotent preradicals is not even weakly idempotent, which also implies that the infimum of weakly idempotent preradicals is not weakly idempotent. The last said that $R$-wid is not a sublattice of $R$-pr and $R$-$\operatorname{cid}$ is not a sublattice of $R$-wid.

As always $R$-wid has a natural way to be described as a complete lattice, that is for any family $\{\sigma_{i}\}_{i\in I}$ of weakly idempotent preradicals the supremum is the usual supremum in $R$-pr, but the infimum results $(\bigwedge_{i\in I}\sigma_{i})^{*}.$The next proposition tells that the operator ${}^{*}$ over $R$-pr is an interior operator.

\begin{proposition}
The assignation ${}^{*}:R$-pr--$\to R$-pr given by $\sigma\mapsto\sigma^{*}$ for any preradical over $R$-Mod $\sigma$ is a monotone, deflatory and idempotent operator over $R$-pr.
\end{proposition}

From the fact that $R$-$\operatorname{id}_{\subseteq}^{-}R$-$\operatorname{cid}$ follows that $\sigma^{\circ}(M)\leq\sigma^{*}(M)$ for any left $R$-module $M$. Then

\begin{remark}
Let $\sigma$ be an essentially idempotent preradical over $R$-Mod and let $M$ be a left $R$-module. If $\sigma(M)=M$ then $\sigma^{*}(M)=M$.
\end{remark}

\begin{remark}
Let $\sigma$ be a preradical over $R$-Mod. Then $\mathbb{T}_{\sigma}=\mathbb{T}_{\sigma^{*}}=\mathbb{T}_{\widehat{\sigma}}$.
\end{remark}

\begin{remark}
Let $\sigma$ be a preradical over $R$-Mod. Then $\widehat{\sigma^{*}}=\widehat{\sigma}$.
\end{remark}

\begin{proposition}
Let $\sigma$ be a weakly idempotent preradical over $R$-Mod. Then $\bar{\sigma}$ is a weakly idempotent preradical.
\end{proposition}

\textbf{Proof.} If $\bar{\sigma}\neq 0$, then $\sigma\neq 0$ it follows that $\widehat{\sigma}\neq 0$ and as $\widehat{\sigma}\leq(\widehat{\bar{\sigma}})$ the desired result is obtained.

\begin{proposition}
Let $S$ be a simple left $R$-module. If $\alpha_{S}^{E(S)}$ is weakly idempotent then $\alpha_{S}^{E(S)}$ is idempotent.
\end{proposition}

\begin{corollary}
Let $R$ be a ring. Then $R$ is a V-ring if and only if every atom in $R$-pr is weakly idempotent.
\end{corollary}

\begin{remark}
Let $\{\sigma_{i}\}_{i\in I}$ be a family of preradicals over $R$-Mod such that $\bigwedge_{i\in I}\sigma$ is weakly idempotent. Then $\bigwedge_{i\in I}\widehat{\sigma}$ is weakly idempotent.
\end{remark}

The complement of $R$-wid in $R$-pr is the class of all no zero preradicals such that the idempotent associated is the zero preradical, a preradical with this property is called a strongly nilpotent preradical. The class of all strongly nilpotent preradicals is denoted by $R$-stn. It is remarked that $R$-stn is closed under infimum but as not for all preradical $\sigma$ there is not necessarily a strongly nilpotent preradical $\tau$ such that $\sigma\leq\tau$ then there is not possible to construct the least strongly nilpotent preradical over $\sigma$. Easily all atoms are idempotents or strongly nilpotents. Also the class $R$-stn is closed under subpreradicals, that is, let $\sigma$ and $\tau$ be preradicals over $R$-Mod such that $\sigma\leq\tau$ and $\tau\in R$-stn then $\sigma\in R$-stn.

\section{Prehereditary Preradicals}
Let $\sigma$ be a preradical over $R$-Mod, $\sigma$ is called a prehereditary preradical if the class $\mathbb{T}_{\sigma}$ is hereditary. It is observed that $\sigma$ is prehereditary if and only if $\widehat{\sigma}$ is hereditary, since $\mathbb{T}_{\sigma}=\mathbb{T}_{\widehat{\sigma}}$. It is remembered that a preradical $\sigma$ is hereditary and only if $\sigma$ is idempotent and $\mathbb{T}_{\sigma}$ is hereditary, this means that prehereditary preradicals are a generalization of hereditary preradicals, the requirement of being idempotent is discarded. . As for any family of preradicals $\\{\sigma_{i}\\}_{i\in I}$, $\mathbb{T}_{\bigwedge_{i\in I}\sigma_{i}}=\bigcap_{i\in I}\mathbb{T}_{\sigma_{i}}$, and the infimum of hereditary pretorsion classes is an hereditary pretorsion class, it follows:

\begin{proposition}
Let $\{\sigma_{i}\}_{i\in I}$ a family of prehereditary preradicals over $R$-Mod. Then $\bigwedge_{i\in I}\sigma_{i}$ is a prehereditary preradical.
\end{proposition}

Let $\sigma$ be a preradical over $R$-Mod. It is possible to construct the least prehereditary preradical over $\sigma$, it will be denoted by $\sigma^{\square}$, and it results the infimum of all prehereditary preradicals over $\sigma$.The class of all prehereditary preradicals is denoted by $R$-pher, and it follows that $R$-lep$\subseteq R$-pher, this implies that $\sigma^{\square}\leq\widetilde{\sigma}$ for any preradical over $R$-Mod $\sigma$. It is observed that a preradical $\sigma$ is prehereditary if and only if $\widetilde{\sigma^{\square}}=\sigma$.The next proposition tells that the operator ${}^{\square}$ over $R$-pr is a closure operator.

\begin{proposition}
The assignation ${}^{\square}:R$-pr$\longrightarrow R$-pr given by $\sigma\mapsto\sigma^{\square}$ for any preradical over $R$-Mod is a monotone, inflatory and idempotent operator over $R$-pr.
\end{proposition}

\begin{remark}
Let $\sigma$ be a preradical over $R$-Mod. Then $\widetilde{\sigma^{\square}}=\widetilde{\sigma}$.
\end{remark}

For a left max ring $R$ (a ring is left max if every no zero left module has a maximal submodule), if it is considered the Jacobson radical $J$ then in general $J$ is not idempotent and ${\mathbb{T}}_{J}=\{0\}$ which means that $J$ is a prehereditary preradical. In particular for any prime $p$ and positive integer $n$, ${\mathbb{Z}}_{p^{n}}$ is a max ring and the Jacobson radical is $\omega_{0}^{p^{n-1}{\mathbb{Z}}_{p^{n}}}$. Moreover $\omega_{0}^{p^{k}{\mathbb{Z}}_{p^{n}}}$ with $k=1,...,n-1$ is a prehereditary preradical which is not an hereditary preradical.

\begin{proposition}
Let $J$ be the Jacobson radical. It is equivalent for $J$:
\begin{enumerate}
\item $J$ is a prehereditary preradical.
\item ${\mathbb{T}}_{J}=\{0\}$.
\item $\widehat{J}=0$.
\item $R$ is a left max ring.
\end{enumerate}
\end{proposition}

\textbf{Proof.} The interesting part is (1) implies (2) and the others are quite obvious. Let $M$ be a no zero left $R$-module and let $x\in M$ a no zero element, as $Rx$ has a maximal left submodule $J(Rx)\neq Rx$ then $J(M)\neq M$.

\begin{lemma}
Let $\sigma$ be a preradical over $R$-Mod. Then $\sigma\leq\widehat{\sigma^{\square}}$ if and only if $\widetilde{\sigma}=\sigma^{\square}$.
\end{lemma}

\textbf{Proof.} If $\sigma\leq\widehat{\sigma^{\square}}$ and as $\widehat{\sigma^{\square}}$ is a left exact preradical, then $\widetilde{\sigma}\leq\widehat{\sigma^{\square}}\leq\sigma^{\square}$ which implies $\widetilde{\sigma}=\sigma^{\square}$. In the case that $\widetilde{\sigma}=\sigma^{\square}$, it follows that $\sigma\leq\sigma^{\square}=\widetilde{\sigma}=\widehat{\widetilde{\sigma}}= \widehat{\sigma^{\square}}$ . $\blacksquare$

The next proposition says that the operator ${}^{\square}$ preserves the idempotency.

\begin{proposition}
Let $\sigma$ be a preradical over $R$-Mod. If $\sigma$ is idempotent then $\sigma^{\square}$ is left exact.
\end{proposition}

\textbf{Proof.} As $\sigma\leq\sigma^{\square}$ then $\sigma\leq\widehat{\sigma}\leq\widehat{\sigma^{\square}}$ and by the previous proposition the result follows. $\blacksquare$

For any preradical $\sigma$ over $R$-Mod it is assigned the set of left ideals ${\mathbb{I}}_{\sigma}=\{_{R}I\leq R\mid R/I\in{\mathbb{T}}_{\sigma}\}$, this set is always closed under over ideals, which means, that if $I\in{\mathbb{I}}_{\sigma}$ and ${}_{R}J\leq R$ is such that $I\subseteq J$ then $J\in{\mathbb{I}}_{\sigma}$, this follows since ${\mathbb{T}}_{\sigma}$ is closed under quotients. The next proposition give sufficient and necessary conditions for the set ${\mathbb{I}}_{\sigma}$ to be a linear filter.

\begin{proposition}
Let $\sigma$ be a preradical over $R$-Mod. Then $\sigma$ is prehereditary if and only if ${\mathbb{I}}_{\sigma}$ is a linear filter.
\end{proposition}

\textbf{Proof.} If $\sigma$ is prehereditary, the proof that ${\mathbb{I}}_{\sigma}$ is a linear filter is the same proof as when $\sigma$ is left exact. Now, let $\tau$ be the left exact preradical induced by ${\mathbb{I}}_{\sigma}$ and let $M$ be a $\tau$-torsion left $R$-module, then $ann(x)\in{\mathbb{I}}_{\sigma}$ for any $x\in M$, therefore $\sigma(Rx)=Rx$ for any $x\in M$, from this $\sigma(M)=M$ and ${\mathbb{T}}_{\tau}\subseteq{\mathbb{T}}_{\sigma}$. Let $M$ a left $R$-module which is not $\tau$-torsion, then there is $x\in M$ with $\tau(Rx)\neq Rx$, as $Rx\cong R/ann(x)$ and ${\mathbb{I}}_{\sigma}={\mathbb{I}}_{\tau}$ it follows that $\sigma(Rx)\neq Rx$ and $\sigma(M)\neq M$, hence $\mathbb{T}_{\tau}\subseteq \mathbb{T}_\sigma$. $\blacksquare$

\begin{corollary}
Let $\sigma$ be a preradical over $R$-Mod. Then ${\mathbb{T}}_{\sigma}$ is an hereditary torsion class if and only if ${\mathbb{I}}_{\sigma}$ is a Gabriel filter.
\end{corollary}

\textbf{Proof.} From the previous result ${\mathbb{I}}_{\sigma}$ is a linear filter and the fact that ${\mathbb{T}}_{\sigma}$ is an hereditary torsion class implies that ${\mathbb{I}}_{\sigma}$ is a Gabriel filter. As in the previous proof , ${\mathbb{T}}_{\tau}\subseteq{\mathbb{T}}_{\sigma}$ where $\tau$ is the left exact radical induced by ${\mathbb{I}}_{\sigma}$. $\blacksquare$

Let $\sigma$ be a preradical over $R$-Mod it is called costable if ${\mathbb{F}}_{\sigma}$ is closed under injective hulls.

\begin{proposition}
Let $\sigma$ be a radical over $R$-Mod. If $\sigma$ is costable then $\sigma$ is left exact.
\end{proposition}

\textbf{Proof.} See [1]. $\blacksquare$

\begin{proposition}
Let $\sigma$ be an essentially idempotent prehereditary preradical over $R$-Mod. Then $\sigma$ is costable.
\end{proposition}

\textbf{Proof.} Let $M$ be a left $R$-module such that $\sigma(M)=0$, then $\widehat{\sigma}(M)=0$, since $\widehat{\sigma}(M)=\widehat{E(M)})\cap M$ and $M\trianglelefteq E(M)$ it follows that $\widehat{\sigma}(E(M))=0$, as $\sigma$ is essentially idempotent $\sigma(EM)=0$.

\begin{proposition}
Let $\{\sigma_{i}\}_{i\in I}$ be a family of prehereditary radicals over $R$-Mod. Then $\bigwedge_{i\in I}\sigma_{i}$ is a prehereditary radical.
\end{proposition}

\begin{proposition}
Let $\sigma$ be a preradical over $R$-Mod and let $M$ be a left $R$-module. If $M\in\mathbb{T}_{p}$ then $M\subseteq\sigma^{\square}(N)$ for any left $R$-module $N$ such that $M\leq N$, in particular $M\subseteq\sigma^{\square}(E(M))$.
\end{proposition}

A counterexample that the inverse proposition is not valid, let $p$ and $q$ be different primes and it is defined $\sigma=\alpha_{Z_{p}}^{Z_{p^{\infty}}}\vee\alpha_{Z_{q}}^{\bar{Z}_{q^{\infty}}}$, as $\sigma$ is prehereditary (since $\widehat{\sigma}=0$) it follows that $\sigma^{\square}=\sigma$, and $\sigma(\mathbb{Z}_{q^{\infty}})=\mathbb{Z}_{q}$ then $\mathbb{Z}_{q}\subseteq\sigma(E(\mathbb{Z}_{q}))$, but $\sigma(\mathbb{Z}_{q})=0$.

\begin{proposition}
Let $S$ be a simple left $R$-module. Then $\alpha_{S}^{S}(E(M))=\alpha_{S}^{S}(M)$ for any left $R$-module $M$.
\end{proposition}

\begin{corollary}
Let $\sigma$ be an atom in $R$-pr. Then $\sigma$ is prehereditary.
\end{corollary}

\textbf{Proof.} As $\sigma$ is an atom it must be of the form $\alpha_{S}^{E(S)}$ for some simple left $R$-module $S$, since $\sigma$ is an atom $\widehat{\sigma}$ must be $0$ or $\sigma$, in the first case $\sigma$ is prehereditary and in the second case $\sigma$ is idempotent which means that $S$ is injective, so by the previous proposition $\sigma(M)=\alpha_{S}^{E(S)}(M)=\alpha_{S}^{S}(M)=\alpha_{S}^{S}(E(M))\cap M= \sigma(E(M))\cap M$ which means that $\sigma$ is left exact, therefore prehereditary. $\blacksquare$

\begin{lemma}
Let $M$ be a left $R$-module and let $N$ be a fully invariant submodule of $M$. If $M$ is $L$-injective for all $L\in\mathbb{T}_{\omega^{M}_{N}}$ then $\omega^{M}_{N}$ is a prehereditary preradical.
\end{lemma}

\textbf{Proof.} Let $K$ be a $\omega_{N}^{M}$-torsion module, let $L$ be a submodule of $K$ and let $f:L\longrightarrow M$ be an $R$-morphism, so there is $g:K\longrightarrow M$ such that $g|_{L}=f$ which means that $f^{-1}(N)=g^{-1}(N)\cap L=K\cap L=L$ implying that $L$ is $\omega_{N}^{M}$-torsion. $\blacksquare$

\begin{proposition}
Let $\sigma$ be a left exact preradical over $R$-Mod and let $\tau$ be a strongly nilpotent preradical over $R$-Mod. If $\sigma\wedge\tau=0$ then $\sigma\vee\tau$ is a prehereditary preradical.
\end{proposition}

\textbf{Proof.} Is is noticed that $\widehat{\sigma\vee\tau}=\widehat{\sigma}\vee\widehat{\tau}=\sigma$. $\blacksquare$

The last result tells how to construct a infinite family of non trivial prehereditary preradicals (understanding by non trivial as no left exact), let $p$ and $q$ be different primes then by the previous proposition $\alpha_{Zp}^{Zp}\vee\alpha_{Zq}^{Zq^{\infty}}$ is a prehereditary predical which is not idempotent (meaning that is not left exact), neither is radical, and its torsion class is non trivial.

\section{Essentiality with respect to a Preradical}

\begin{proposition}
Let $M$ be a non singular left $R$-module and let $N$ a submodule of $M$. Then $N\unlhd M$ if and only if $M/N$ is singular.
\end{proposition}

\begin{proposition}
Let $M$ be a left $R$-module and let $N$ and $K$ be submodules of $M$. Then:
\begin{enumerate}
\item Let $x\in M$. If $N\unlhd M$ then $(N:x)\unlhd R$.
\item If $K\unlhd M$ and $L\unlhd M$ then $N\cap K\unlhd M$.
\item If $K\unlhd N$ and $N\unlhd M$ then $K\unlhd M$.
\item If $K\unlhd M$ and $K\leq N$ then $K\unlhd N$ and $N\unlhd M$.
\end{enumerate}
\end{proposition}

With the previous two propositions it is possible to think in a kind of essentiality respect to a preradical $\sigma$, this let $M$ be a left $R$-module and $N$ a submodule of $M$, it is said that $N$ is $\sigma$-dense in $M$ if $\sigma(M/N)=M/N$, this fact is denoted by $N\unlhd_{\sigma}M$, it may be thought as $\sigma$-essentiality, in fact if $\sigma$ is prehereditary most of the last properties are preserved.

\begin{proposition}
Let $M$ be a left $R$-module, let $N$ and $K$ be submodules of $M$ and let $\sigma$ be a preradical over $R$-Mod . Then:
\begin{enumerate}
\item If $K\unlhd_{\sigma}M$ and $K\leq N$ then $N\unlhd_{\sigma}M$. When $\sigma$ is prehereditary,
\item Let $x\in M$. If $N\unlhd_{\sigma}M$ then $(N:x)\unlhd_{\sigma}R$.
\item If $K\unlhd_{\sigma}M$ and $L\unlhd\sigma M$ then $N\cap K\unlhd_{\sigma}M$.
\item If $K\unlhd_{\sigma}M$ and $K\leq N$ then $K\unlhd_{\sigma}N$.
\item If $N\leq M$ and $K\unlhd_{\sigma}M$ then $K\cap N\unlhd_{\sigma}N$. When $\sigma$ is essentially coidempotent,
\item If $K\unlhd_{\sigma}N$ and $N\unlhd_{\sigma}M$ then $K\unlhd_{\sigma}M$.
\end{enumerate}
\end{proposition}

\textbf{Proof.}
\begin{enumerate}
\item Follows from the fact that ${\mathbb{T}}_{\sigma}$ is closed under quotients.
\item It is considered the next equalities $R/(N:x)=R/ann(x+N)\cong R(x+N)\leq M/N$.
\item Let $\pi:M\longrightarrow M/N\times M/K$ be the morphism induced by the canonical projections, as $M/N\times M/K$ is a $\sigma$-torsion left $R$-module, $\ker\pi=N\cap K$ and $M/(N\cap K)$ is isomorphic to a submodule of $M/N\times M/K$ then $M/(N\cap K)$ is a $\sigma$-torsion left $R$-module.
\item Follows from the fact that ${\mathbb{T}}_{\sigma}$ is closed under submodules.
\item Follows by the second isomorphism theorem and the fact that ${\mathbb{T}}_{\sigma}$ is closed under submodules.
\item Follows from the fact that ${\mathbb{T}}_{\sigma}$ is closed under extensions.
\end{enumerate}

\begin{proposition}
Let $\sigma$ and $\tau$ be preradicals over $R$-Mod, let $M$ be a left $R$-module and let $N$ be a submodule of $M$. If $\sigma\leq\tau$ and $N\unlhd_{\sigma}M$ then $N\unlhd_{\tau}M$.
\end{proposition}

\textbf{Proof.} Follows from the fact that ${\mathbb{T}}_{\sigma}\subseteq{\mathbb{T}}_{\tau}$.

\begin{proposition}
Let $\sigma$ be a preradical over $R$-Mod, let $M$ be a left $R$-module and let $N$ be a submodule of $M$. Then $N\unlhd_{\sigma}M$ if and only if $N\unlhd_{\widehat{\sigma}}M$.
\end{proposition}

\textbf{Proof.} Follows from the fact that $\mathbb{T}_{\sigma}=\mathbb{T}_{\widetilde{\sigma}}$. $\blacksquare$

\begin{proposition}
Let $\sigma$ be a prehereditary preradical over $R$-Mod, let $M$ be a left $R$-module and let $N$ be a submodule of $M$. If $M$ is $\sigma$-torsionfree and $N\unlhd_{\sigma}M$ then $N\unlhd M$.
\end{proposition}

\textbf{Proof.} Let $x\in M$, then by the second isomorphism theorem $Rx/(Rx\cap N)\cong(Rx+N)/N$ which implies that $Rx+N\unlhd_{\sigma}M$ and $Rx/(Rx\cap N)$ is of $\sigma$-torsion, since $Rx$ is $\sigma$-torsionfree if $Rx\cap N=0$ then $Rx$ is $\sigma$-torsion and this implies $Rx=0$ and $x=0$, which means that if $x\neq 0$ then $Rx\cap N\neq 0$. $\blacksquare$
\section{Pure Submodules respect to a Preradical}
Let $\sigma$ be a preradical over $R$-Mod, let $M$ be a left $R$-module and let $N$ be a submodule of $M$, it is said that $N$ is $\sigma$-pure submodule of $M$ if $M/N$ is $\sigma$-torsionfree.

\begin{proposition}
Let $\sigma$ be a preradical over $R$-Mod, let $M$ be a left $R$-module and let $\{M_{i}\}_{i\in I}$ be a family of $\sigma$-pure submodules of $M$. Then $\bigcap_{i\in I}M_{i}$ is a $\sigma$-pure submodule of $M$.
\end{proposition}

\textbf{Proof.} Let $\pi:M\longrightarrow\prod_{i\in I}M/M_{i}$ be the morphism induced by canonical projections, since $\prod_{i\in I}M/M_{i}$ is $\sigma$-torsionfree in follows that $M/\bigcap_{i\in I}M_{i}$ is $\sigma$-torsionfree. $\blacksquare$

For a submodule $N$ of a left $R$-module $M$, it should be considered the least $\sigma$-pure submodule of $M$ that contains $N$, it is denoted by $N^{M}_{\sigma}$, it is described by

\[N^{M}_{\sigma}=\bigcap\{K\leq M\mid N\leq K,M/K\in\mathbb{F}_{\sigma}\}\]

as the last proposition states it is $\sigma$-pure in $M$ and contains $N$. The submodule $N^{M}_{\sigma}$ is called the $\sigma$-purification of $N$ in $M$.

\begin{proposition}
Let $\sigma$ be a preradical over $R$-Mod, let $M$ be a left $R$-module and let $N$ be a submodule of $M$. Then $\bar{\sigma}(M/N)=N^{M}_{\sigma}/N$.
\end{proposition}

\begin{corollary}
Let $\sigma$ be a preradical over $R$-Mod, let $M$ be a left $R$-module and let $N$ be a submodule of $M$. Then $N^{M}_{\sigma}=N^{M}_{\widetilde{\sigma}}$.
\end{corollary}

\begin{proposition}
Let $\sigma$ be a preradical over $R$-Mod, let $M$ be a left $R$-module and let $N$ be a submodule of $M$. Then $N$ is $\sigma$-pure in $M$ if and only if $N=N^{M}_{\sigma}$.
\end{proposition}

\textbf{Proof.} As $\sigma(M/N)=0$ implies $\bar{\sigma}(M/N)=0$ the result follows. $\blacksquare$

\begin{remark}
Let $\sigma$ be a preradical over $R$-Mod and let $M$ be a left $R$-module. Then $\bar{\sigma}(M)=0^{M}_{\sigma}$.
\end{remark}

\begin{lemma}
Let $\sigma$ be a preradical over $R$-Mod, let $M$ be a left $R$-module and let $N$ be a submodule of $M$. If $N$ is $\sigma$-pure in $M$ and $N\in\mathbb{T}_{\sigma}$ then $\sigma(M)=N$.
\end{lemma}

\textbf{Proof.} First as $N$ is $\sigma$-pure in $M$ then $\sigma(M/N)=0$, this way

\[\sigma(M) \subseteq\bar{\sigma}(M)\] \[=\bigcap\{K\leq M\mid\sigma(M/K)=0\}\subseteq N\]

By other side $N\leq M$ implies $N=\sigma(N)\leq\sigma(M)$. $\blacksquare$

\begin{remark}
Let $\sigma$ be an idempotent radical over $R$-Mod, let $M$ be a left $R$-module and let $N$ be a submodule of $M$. Then $N^{M}_{\sigma}/N$ is a $\sigma$-torsion module.
\end{remark}

\section{Injectivity respect to a Preradical}
Let $\sigma$ be a preradical over $R$-Mod and let $M$ be a left $R$-module, it is said that $M$ is $\sigma$-injective if $f:K\longrightarrow M$ is a $R$-morphism and $K\trianglelefteq_{\sigma}N$ then there is a morphism $g:N\longrightarrow M$ with $g|_{K}=f$. This concept is a generalization of injectivity respect an hereditary torsion theory, as a reference is the book (2) .The first thing that is observed is that if $M$ is also a $\sigma$-torsion module then $M$ is quasi-injective.

\begin{proposition}
Let $\sigma$ be a preradical over $R$-Mod and let $M$ be a left $R$-module. Then (1) and (2) are equivalent and imply (3), (3) implies (4) and (5), (4) implies (6) and (5) implies (6). If $\sigma$ is an idempotent radical then (1), (2), (3) and (5) are equivalent. If $\sigma$ is prehereditary then (5) and (6) are equivalent. If $\sigma$ is a left exact radical then all are equivalent.
\begin{enumerate}
\item $M$ is $\sigma$-pure in $E(M)$
\item If $M$ is a submodule of a left $R$-module $N$, then there exist a $\sigma$-pure submodule $K$ of $N$ that contains $M$ and $M$ is a direct summand of $K$.
\item $Ext^{1}_{R}(N,M)=0$ for any $\sigma$-torsion left $R$-module $N$.
\item $Ext^{1}_{R}(R/I,M)=0$ for any $I\in\mathbb{I}_{\sigma}$.
\item $M$ is $\sigma$-injective
\item $M$ is $\sigma$-injective respect to $R$
\end{enumerate}
\end{proposition}

\textbf{Proof.}
(1) $\Rightarrow$ (2) If $M\leq N$ then $E(N)=E(M)\oplus L$. It is put $L^{\prime}=N\cap L$ then $N/(M\oplus L^{\prime})$ is isomorphic to a submodule of $E(N)/(M\oplus L)$, since it is considered the morphism $f:N\longrightarrow E(N)/(M\oplus L)$ with $f=gh$ where $g:E(N)\longrightarrow E(N)/(M\oplus L)$ is the canonical projection and $h:N\longrightarrow E(N)$ is the canonical inclusion, so ker $f=N\cap(M\oplus L)=M\oplus L^{\prime}$. By other side $E(N)/(M\oplus L)\cong(E(M)\oplus L)/(M\oplus L)\cong E(M)/M$ which by hypothesis is $\sigma$-torsionfree, that is why $N/(M\oplus L^{\prime})$ is $\sigma$-torsionfree and $K=M\oplus L^{\prime}$ is $\sigma$-pure in $N$.

(2) $\Rightarrow$ (1) As $M\leq E(M)$, by hypothesis then there is $K\leq E(M)$ with $M\oplus K$$\sigma$-pure in $E(M)$, but $M\trianglelefteq E(M)$ which means that $K=0$, therefore $M$ is $\sigma$-pure in $E(M)$.

(1) $\Rightarrow$ (3) It is considered the short exact sequence

\[0\longrightarrow M\longrightarrow E(M)\longrightarrow E(M)/M\longrightarrow 0\]

and it is obtained a exact sequence

\[Hom_{R}(N,E(M)/M)\longrightarrow Ext^{1}_{R}(N,M)\longrightarrow Ext^{1}_{R}(N ,E(M))\]

where $N$ is a $\sigma$-torsion left $R$-module, as $E(M)/M$ is $\sigma$-torsionfree, which implies $Hom_{R}(N,E(M)/M)=0$ and as $E(M)$ is injective $Ext^{1}_{R}(N,E(M))=0$, from this follows that $Ext^{1}_{R}(N,M)=0$.

(3) $\Rightarrow$ (5) It is taken a short exact sequence $0\longrightarrow N^{\prime}\longrightarrow N\longrightarrow N/N^{\prime} \longrightarrow 0$ such that $N/N^{\prime}$ is a $\sigma$-torsion module, it is induced the following short exact sequence

\[0\longrightarrow Hom_{R}(N/N^{\prime},M)\longrightarrow Hom_{R}(N,M) \longrightarrow Hom_{R}(N,M)\longrightarrow Ext^{1}_{R}(N/N^{\prime},M)\]

and by hypothesis the last module is zero.

(2) $\Rightarrow$ (4),(4) $\Rightarrow$ (6) and (5) $\Rightarrow$ (6) are obvious.

(6) $\Rightarrow$ (5) In the same way as the proof of the Baer's criterion.

(5) $\Rightarrow$ (1) As $M\unlhd_{\sigma}M^{E(M)}_{\sigma}$ since $\sigma$ is an idempotent radical, there is an $R$-morphism $\alpha:M^{E(M)}_{\sigma}\longrightarrow M$ such that $\alpha|M=1_{M}$ which implies that $\alpha$ is an epimorphism, it is noticed that $\ker\alpha\cap=\ker 1_{M}=0$ since $M\unlhd M^{E(M)}_{\sigma}$ it follows that $\alpha$ is a monomorphism, so $M=M^{E(M)}_{\sigma}$.

By the last proposition it is observed that it is sufficient to ask to a preradical to be an idempotent radical to speak about relative injectivity the only thing that may not be assured is the Baer's criterion and in the other hand it is sufficient to ask a preradical to be prehereditary to have Baer's criterion.

Let $\sigma$ be a preradical over $R$-Mod and let $M$ be a left $R$-module, it is defined the $\sigma$-injective hull of $M$ as $M^{E(M)}_{\sigma}$ and it is denoted by $E_{\sigma}(M)$.

\begin{remark}
Let $\sigma$ be a preradical over $R$-Mod and let $M$ be a left $R$-module. Then:
\begin{enumerate}
\item $E_{\sigma}(M)$ is $\sigma$-injective
\item $M\unlhd E_{\sigma}(M)$
\item If $\sigma$ is an idempotent radical then $M\unlhd_{\sigma}E_{\sigma}(M)$.
\end{enumerate}
\end{remark}

This three properties characterizes the $\sigma$-injective hull as the next proposition tells. It is an analogous characterization of the usual injective hull as an injective essential extension, but now it is asked to be $\sigma$-dense extension and the hypothesis over $\sigma$ to be an idempotent radical.

\begin{proposition}
Let $\sigma$ be a preradical over $R$-Mod, let $K$ be a $\sigma$-injective left $R$-module and let $M$ be a $\sigma$-dense dense submodule of $K$. If $\sigma$ is an idempotent radical then $K=E_{\sigma}(M)$.
\end{proposition}

\textbf{Proof.} As $M$ is essential in $K$ without loss of generality it may reduced to the case when $K\leq E(M)$, so $E(K)=E(M)$ and $K$ is $\sigma$-pure in $E(M)$ then by lemma 3 the result is followed.

\begin{remark}
Let $\sigma$ be an idempotent radical over $R$-Mod and let $M$ be a left $R$-module. Then $M$ is $\sigma$-injective if and only if $E_{\sigma}(M)=M$.
\end{remark}

\begin{proposition}
Let $\sigma$ be an idempotent radical over $R$-Mod, let $M$ be a left $R$-module and let $N$ be a submodule of $M$. If $M$ is $\sigma$-injective and $N$ is $\sigma$-pure in $M$ then $N$ is $\sigma$-injective.
\end{proposition}

\textbf{Proof.} Let $K$ be a $\sigma$-torsion left $R$-module and it is considered the next short exact sequence:

\[Hom_{R}(K,M/N)\longrightarrow Ext^{1}_{R}(K,N)\longrightarrow Ext^{1}_{R}(K,M)\]

As $Ext^{1}_{R}(K,M)=0$ and $Hom_{R}(K,M/N)=0$ since $M$ is $\sigma$-injective and $M/N$ is $\sigma$-torsion free it follows that $Ext^{1}_{R}(K,N)=0$, therefore $N$ is $\sigma$-injective. $\blacksquare$

\begin{remark}
Let $\sigma$ be an idempotent radical over $R$-Mod and let $M$ be a left $R$-module. If $M$ is $\sigma$-injective $\sigma$-torsion module the $\sigma(E(M))=M$.
\end{remark}

\begin{remark}
Let $\sigma$ be a radical over $R$-Mod and let $M$ be a left $R$-module. If $\sigma(E(M))=M$ then $M$ is $\sigma$-injective.
\end{remark}

Let $M$ a left $R$-module, it is defined $\Omega(M)$ as the set of all left ideals that contain $ann(x)$ for some $x\in M$.

\begin{lemma}[Technical]
Let $M$ be a left $R$-module. Then $M$ is quasinductive if and only if for any left ideal $L$ and for any $R$-morphism $\alpha:L\longrightarrow M$ with $\ker\alpha\in\Omega(M)$ there is an $R$-morphism $\beta:R\longrightarrow M$ such that $\beta|_{L}=\alpha$.
\end{lemma}

\textbf{Proof.} [3, lemma 2]. $\blacksquare$

The previous lemma is used in the proof of the following proposition.

\begin{proposition}
Let $\sigma$ be a preradical over $R$-Mod and let $\{M_{i}\}_{i\in I}$ be a family of left $R$-modules. Then $\prod_{i\in I}M_{i}$ is $\sigma$-injective if and only if $M_{i}$ is $\sigma$-injective for any $i\in I$.
\end{proposition}

\begin{proposition}
Let $\sigma$ be a preradical over $R$-Mod and let $M$ be a $\sigma$-torsion left $R$-module. If $\sigma$ is an idempotent radical the (1) implies (2), if $\sigma$ is prehereditary then (2) implies (1) and if $\sigma$ is a left exact radical (1) and (2) are equivalent.
\begin{enumerate}
\item $M$ is $\sigma$-injective
\item \begin{enumerate}
\item $M$ is quasi-injective.
\item If $I\in\mathbb{I}_{\sigma}$ and $I^{\prime}$ is a left ideal such that $I^{\prime}\subseteq I$ and $I/I^{\prime}$ can be embedded in $M$ then $I^{\prime}=I\cap ann(x))$ for some $x\in M$.
\end{enumerate}
\end{enumerate}
\end{proposition}

\textbf{Proof.} The arguments of the proposition (4.2) of (4). $\blacksquare$

\begin{proposition}
Let $M$ be a quasi-injective left $R$-module. If $\omega_{M}^{E(M)}$ is a radical then $M$ is $\sigma$-injective for any preradical $\sigma$ such that $\sigma(E(M))=M$.
\end{proposition}

\textbf{Proof.} If $\omega_{M}^{E(M)}$ is a radical, as $\sigma(E(M))=M$ implies $\sigma\leq\omega_{M}^{E(M)}$, so $\sigma(E(M)/M)\leq\omega_{M}^{E(M)}(E(M)/M)=0$ which means that $M$ is $\sigma$-pure in $E(M)$, therefore $\sigma$-injective. $\blacksquare$

\textbf{Examples.} Let $\sigma$ be a preradical over $R$-Mod and let $M$ be a left $R$-module, as it is seen $E_{\sigma}(M)/M=\bar{\sigma}(E(M)/M)$. Let $R$ be the ring of the integers $\mathbb{Z}$ and it is considered $\sigma=Soc,t,d,J$ where $Soc$ is the socle, $t$ the torsion part, $d$ the divisible part and $J$ the Jacobson radical then $E_{\sigma}(\mathbb{Z})=\mathbb{Q}$ and $E_{\sigma}(\mathbb{Z}_{p^{k}})=\mathbb{Z}_{p^{\infty}}$ with $p$ a prime number and $k$ a natural number. But if $\sigma=\alpha\frac{\mathbb{Z}_{p}}{\mathbb{Z}_{p}}$ with $p$ a prime number then $E_{\sigma}(\mathbb{Z})=\{\frac{a}{p^{m}}\in\mathbb{Q}\mid a,m\in\mathbb{N}\}$, $E_{\sigma}(\mathbb{Z}_{q^{k}})=\mathbb{Z}_{q^{\infty}}$ when $p=q$ and $E_{\sigma}(\mathbb{Z}_{q^{k}})=\mathbb{Z}_{q^{k}}$ if $p\neq q$ with $q$ a prime number.

\section{Pseudocomplemented Submodules relative to a Preradical}
Let $\sigma$ be a preradical over $R$-Mod, let $M$ be a left $R$-module and let $N$ be a submodule of $M$, it is said that $N$ is $\sigma$-pseudocomplemented in $M$ if there is a submodule $K$ of $M$ such that $N\cap K=0$, $N\oplus K\unlhd M$ and $N\oplus K$ is $\sigma$-dense in $M$ , the submodule $K$ is called a $\sigma$-pseudocomplement of $N$ in $M$.

\begin{proposition}
Let $\sigma$ be a essentially coidempotent preradical over $R$-Mod and let $M,N$ and $K$ be a left $R$-modules with $K\leq N\leq M$. If $K$ is $\sigma$-pseudocomplemented in $N$ and $N$ is $\sigma$-pseudocomplemented in $M$ then $K$ is $\sigma$-pseudocomplemented in $M$.
\end{proposition}

\textbf{Proof.} By hypothesis there are $K^{\prime}$ submodule of $N$ and $N^{\prime}$ submodule of $M$ such that $K\cap K^{\prime}=0$, $N\cap N^{\prime}=0$, $K\oplus K^{\prime}\trianglelefteq N$, $N\oplus N^{\prime}\trianglelefteq M$, $K\oplus K^{\prime}\trianglelefteq_{\sigma}N$ and $N\oplus N^{\prime}\trianglelefteq_{\sigma}M$. It is proposed $K^{\prime}\oplus N^{\prime}$ as the $\sigma$-pseudocomplement of $K$ in $M$. Immediately $K\oplus K\oplus N^{\prime}\trianglelefteq M$. Next it is considered the following short exact sequence:

\[0\longrightarrow(N\oplus N^{\prime})/(K\oplus K\oplus N^{\prime})\longrightarrow M /(K\oplus K\oplus N^{\prime})\longrightarrow M/(N\oplus N^{\prime})\longrightarrow 0\]

As $(N\oplus N^{\prime})/(K\oplus K\oplus N^{\prime})\cong N/(K\oplus K^{\prime})\in \mathbb{T}_{\sigma}$, $M/(N\oplus N^{\prime})\in\mathbb{T}_{\sigma}$ and $\mathbb{T}_{\sigma}$ is closed under extensions then $M/(K\oplus K\oplus N^{\prime})\in\mathbb{T}_{\sigma}$ and the proposition it is obtained. $\blacksquare$

\begin{proposition}
Let $\sigma$ be a prehereditary preradical over $R$-Mod and let $M,N$ and $K$ be a left $R$-modules with $K\leq N\leq M$. If $K$ is $\sigma$-pseudocomplemented in $M$ then $K$ is $\sigma$-pseudocomplemented in $N$.
\end{proposition}

\textbf{Proof.} By hypothesis there is $K^{\prime}$ submodule of $M$ such that $K\cap K^{\prime}=0$, $K\oplus K^{\prime}\trianglelefteq M$ and $K\oplus K^{\prime}\trianglelefteq_{\sigma}M$. It is proposed $K^{\prime\prime}=N\cap K^{\prime}$ as the $\sigma$-pseudocomplement of $K$ in $N$. First it is obvious that $K^{\prime\prime}\cap N=0$. Second

\[K\oplus K^{\prime\prime} =K\oplus(N\cap K^{\prime})\] \[=N\cap(K\oplus K^{\prime})\trianglelefteq N\]

At last, as $K\oplus K^{\prime}\trianglelefteq_{\sigma}M$ then $K\oplus K^{\prime\prime}=N\cap(K\oplus K^{\prime})\trianglelefteq_{\sigma}N$ as it is desired. $\blacksquare$

Let $\sigma$ be a prehereditary preradical over $R$-Mod and let $M$ be a left $R$-module, $Subp_{\sigma}(M)$ denotes the set of all submodules of $M$ that are $\sigma$-pseudocomplemented.

\begin{remark}
Let $\sigma$ and $\tau$ be preradicals over $R$-Mod and let $N$ and $M$ be left $R$-modules. Then:
\begin{enumerate}
\item $M\in Subp_{\sigma}(M)$.
\item $0\in Subp_{\sigma}(M)$.
\item If $N\trianglelefteq_{\sigma}M$ then $N\in Subp_{\sigma}(M)$.
\item If $M$ is a $\sigma$-torsion module then $Subp_{\sigma}(M)=Sub(M)$.
\item If $\sigma\leq\tau$ then $Subp_{\sigma}(M)\subseteq Subp_{\tau}(M)$.
\end{enumerate}
\end{remark}

Let $\sigma$ be a prehereditary preradical over $R$-Mod, ${\mathbb{E}}_{\sigma}$ denotes the class of all $\sigma$-injective left modules.

\begin{proposition}
Let $\sigma$ and $\tau$ be preradicals over $R$-mod. If $Subp_{\sigma}(M)=Subp_{\tau}(M)$ for any left $R$-module $M$ then ${\mathbb{E}}_{\sigma}={\mathbb{E}}_{\tau}$.
\end{proposition}

\textbf{Proof.} Let $E$ be a $\sigma$-injective left $R$-module, let $M$ be a left $R$-module, $N$ a $\tau$-dense submodule of $M$ and $\alpha:N\longrightarrow E$ an $R$-morphism. First it is observed that $N\in Subp_{\tau}$ then it has a $\sigma$-pseudocomplented $N^{\prime}$ in $M$, then it is considered the morphism $\alpha\oplus 0:N\oplus N^{\prime}\longrightarrow E$ since $N\oplus N^{\prime}\trianglelefteq_{\sigma}M$ then there is a morphism $\beta:M\longrightarrow E$ such that $\beta|_{N\oplus N^{\prime}}=\alpha\oplus 0$, so $\beta|_{N}=\alpha$ which proves that $E$ is $\tau$-injective.

\begin{corollary}
Let $Z$ be the singular preradical and let $E$ be a left $R$-module. Then $E$ is injective if and only if $E$ is $Z$-injective.
\end{corollary}

\section{Absolute $\sigma$-Pure}
Let $\sigma$ be a preradical over $R$-Mod and let $M$ be a left $R$-module. It is said that $M$ is absolutely $\sigma$-pure if $M$ is $\sigma$-torsionfree and $\sigma$-injective.

\begin{proposition}
Let $\sigma$ be a preradical over $R$-Mod and let $M$ be a left $R$-module. Then (1)$\Rightarrow$(2) and if $\sigma$ is idempotent then are equivalents.
\begin{enumerate}
\item $M$ is absolutely $\sigma$-pure.
\item For any left $R$-module $N$, for any $\sigma$-dense submodule of $N$, $K$, and for any $R$-morphism $\alpha:K\longrightarrow M$ there is a unique $R$-morphism $\beta:N\longrightarrow M$ such that $\beta|K=\alpha$.
\end{enumerate}
\end{proposition}

\textbf{Proof.} (1)$\Rightarrow$(2) Let $\beta$ and $\beta^{\prime}$ be $R$-morphisms such that $\beta|K=\alpha$ and $\beta^{\prime}|K=\alpha$. Then $K\leq\ker(\beta-\beta^{\prime})$ so there is a morphism $\gamma:N/K\longrightarrow M$ given by $\gamma(x+K)=(\beta-\beta^{\prime})(x)$ for any $x+K\in N/K$. As $N/K$ is a $\sigma$-torsion module and $M$ is a $\sigma$-torsion free module it follows that $\gamma=0$ and so $\beta=\beta^{\prime}$.

(2)$\Rightarrow$(1) It must seen that $\sigma(M)=0$, so $0$ is $\sigma$-essential submodule of $\sigma(M)$ and there are two morphisms that extend the morphism $0:0\longrightarrow M$, the inclusion $i:\sigma(M)\longrightarrow M$ and $0:\sigma(M)\longrightarrow M$, by the uniqueness $i=0$ and it follows that $\sigma(M)=0$.

\begin{proposition}
Let $\sigma$ be a preradical costable over $R$-Mod and let $M$ be a left $R$-module. If $M$ is $\sigma$-torsion free and $M$ is $\sigma$-pure in any $\sigma$-torsion free module that contains it then $M$ is absolutely $\sigma$-pure.
\end{proposition}

\textbf{Proof.} As $M$ is $\sigma$-torsion free the $E(M)$ is $\sigma$-torsion free, so $M$ is $\sigma$-pure in $E(M)$ which implies that $M$ is $\sigma$-injective.

\begin{proposition}
Let $\sigma$ be an essentially idempotent preradical over $R$-Mod and let $M$ be a left $R$-module. If $M$ is absolutely $\sigma$-pure then $M$ is $\sigma$-torsion free and $M$ is $\sigma$-pure in any $\sigma$-torsion free module that contains it.
\end{proposition}

\textbf{Proof.} Let $M^{\prime}$ be a $\sigma$-torsion free $R$-module that contains $M$, then there is a submodule of $M^{\prime}$, $N$, such that $M\oplus N$ is $\sigma$-pure in $M^{\prime}$. So it is observed the following short exact sequence:

\[0\longrightarrow(M\oplus N)/M\longrightarrow M^{\prime}/M \longrightarrow M^{\prime}/(M\oplus N)\longrightarrow 0\]

Now, as $(M\oplus N)/M\cong N$ which is $\sigma$-torsion free and $M^{\prime}/(M\oplus N)$ is $\sigma$-torsion free then $M^{\prime}/M$ is $\sigma$-torsion free.
\section{Autocostable Preradicals}
Let $\sigma$ be a preradical over $R$-Mod. It is said the $\sigma$ is autocostable if $\mathbb{F}_{\sigma}$ is closed under $\sigma$-injective hulls.

\begin{remark}
Let $\sigma$ be a preradical over $R$-Mod. If $\sigma$ is costable then it is autocostable.
\end{remark}

\begin{proposition}
Let $\sigma$ be a preradical over $R$-Mod. If $\sigma$ is an autocostable essentially idempotent preradical then $\sigma$ is a costable preradical.
\end{proposition}

\textbf{Proof.} Let $M$ be a $\sigma$-torsion free left $R$-module and it is considered the following exact sequence:

\[0\longrightarrow E_{\sigma}(M)\longrightarrow E(M) \longrightarrow E(M)/E_{\sigma}(M)\longrightarrow 0\]

as $E(M)/E_{\sigma}(M)$ is $\sigma$-torsion free then $E(M)$ is $\sigma$-torsion free.

\begin{corollary}
Let $\sigma$ be a preradical over $R$-Mod. If $\sigma$ is an autocostable essentially idempotent radical then $\sigma$ is left exact radical.
\end{corollary}
\section{Localization}
Let $\sigma$ be a preradical over $R$-Mod. It is defined an assignation $Q_{\sigma}$ from $R$-Mod to $R$-Mod as $Q_{\sigma}(M)=E_{\sigma}(M/\sigma(M))$ for any left $R$-module $M$, so it is defined $\eta_{M}^{\sigma}:M\longrightarrow Q_{\sigma}(M)$ as the canonical projection composed with the canonical inclusion. It is observed that if $\sigma$ is a left exact radical then $Q_{\sigma}(M)$ is an absolutely $\sigma$-pure module for any left $R$-module $M$, if $\alpha:M\longrightarrow N$ is an $R$-morphism it induces an $R$-morphism $\tilde{\alpha}:M/\sigma(M)\longrightarrow N/\sigma(N)$ so it is composed with the inclusion of $N/\sigma(N)$ in $Q_{\sigma}(N)$ and as $Q_{\sigma}(N)$ is absolute $\sigma$-pure then there is a unique $R$-morphism $\gamma:Q_{\sigma}(M)\longrightarrow$ such that extends the composition metioned, if it is put $Q_{\sigma}(f)=\gamma$ is straight to check that in this case this assigment makes to $Q_{\sigma}$ an endofunctor over $R$-Mod, the endofunctor is called the localization respect $\sigma$ and has been studied a lot, as references are (4), (5) and (9).

\begin{proposition}
Let $\sigma$ be a left exact radical over $R$-Mod. Then $Q_{\sigma}$ is idempotent and left exact.
\end{proposition}

\begin{proposition}
Let $\sigma$ be a left exact radical over $R$-Mod. Then $\eta^{\sigma}:\operatorname{1}_{R\text{-Mod}}\longrightarrow Q_{\sigma}$ is a natural transformation.
\end{proposition}

\begin{proposition}
Let $\sigma$ be a left exact radical over $R$-Mod and let $M$ be a left $R$-module. Then $\ker\eta_{M}^{\sigma}$ is a $\sigma$-torsion module and $\operatorname{coker}_{M}^{\eta}$ is a $\sigma$-torsion free module.
\end{proposition}

\begin{proposition}
Let $\sigma$ be a left exact radical over $R$-Mod. Then $\eta^{\sigma}\circ Q_{\sigma}=Q_{\sigma}\circ\eta^{\sigma}$.
\end{proposition}

\textbf{Proof.} Let $M$ be a left $R$-module, it is easy to verify that $\eta_{Q_{\sigma}(M)}^{\sigma}=\operatorname{1}_{Q_{\sigma}(M)}$ and $Q_{\sigma}(\eta_{M}^{\sigma})=\operatorname{1}_{Q_{\sigma}(M)}$.

\begin{proposition}
Let $\sigma$ be a preradical over $R$-Mod. Then $\sigma$ is an idempotent preradical if and only if $Q_{\sigma}\circ\sigma=0$
\end{proposition}

\textbf{Proof.} Let $M$ be a left $R$-module then $(Q_{\sigma}\circ\sigma)(M)=E_{\sigma}(\sigma(M)/\sigma^{2}(M))$. $\blacksquare$

\begin{proposition}
Let $\sigma$ be a preradical over $R$-Mod. If $Q_{\sigma}\circ\sigma=\sigma\circ Q_{\sigma}$ then $\sigma$ is an idempotent autocoestable radical.
\end{proposition}

\textbf{Proof.} Is easy to see that $\sigma$ is idempotent, then by the previous proposition $\sigma\circ Q_{\sigma}=0$ which implies that $\sigma(E_{\sigma}(M/\sigma(M)))=0$ for any left $R$-module $M$, so $\sigma(M/\sigma(M))=0$ which means that $\sigma$ is a radical and this implies ${\mathbb{F}}_{\sigma}=\{M/\sigma(M)\mid M\in R\text{-}{\rm Mod}\}$ therefore the class ${\mathbb{F}}_{\sigma}$ is closed under $\sigma$-injective hulls.

\begin{corollary}
Let $\sigma$ be a preradical over $R$-Mod. Then $Q_{\sigma}\circ\sigma=\sigma\circ Q_{\sigma}$ if and only if $\sigma$ is a left exact radical.
\end{corollary}

\textbf{Proof.} All idempotent autocostable radicals are left exact radicals. $\blacksquare$

\section{Bibliography}

\begin{enumerate}
\item Bican, L.; Kepka, T.; Nemec, P., Rings, Modules and Preradicals, Marcel Dekker, New York, 1982.
\item Crivei, S., Injective modules relative to torsion theories, EFES Publishing House, Cluj-Napoca, 2004.
\item Fuchs, L., On quasi-injective modules, Annali della Scuola Normale Superiore di Pisa, Classe di Scienze 3e serie, tome 23, no.4(1969), p. 541-546.
\item Golan, J., Localization of noncommutative rings, Pure and Applied Mathematics, No. 30. Marcel Dekker, Inc., New York, 1975.
\item Golan, J., Torsion theories, Pitman Monographs and Surveys in Pure and Applied Mathematics, 29. Longman Scientific \& Technical, Harlow; John Wiley \& Sons, Inc., New York, 1986.
\item Raggi, F.; Rios, J.; Rincon, H.; Fernandez-Alonso, R.; Signoret, C., The lattice structure of preradicals, Comm. Algebra 30 (2002), no. 3, 1533-1544.
\item Raggi, F.; Rios, J.; Rincon, H.; Fernandez-Alonso, R.; Signoret, C., The lattice structure of preradicals. II, J. Algebra Appl. 1 (2002) no.2, 201-214.
\item Raggi, F.; Rios, J.; Rincon, H.; Fernandez-Alonso, R.; Signoret, C., The lattice structure of preradicals. III, J. Pure Appl. Algebra 190 (2004), no. 1-3, 251-265.
\item Stenstrom, B., Rings of Quotients: An Introduction to Methods of Ring Theory, Springer-Verlag, Berlin, 1975.
\item Wisbauer, R., Foundations of Module and Ring Theory , Gordon and Breach, Reading-Paris, 1991.
\end{enumerate}

\end{document}